# About Evaluation of Complex Dynamical Systems

**Dmytro Polishchuk, Olexandr Polishchuk**

*Pidstryhach Institute for Applied Problems of Mechanics and Mathematics National Academy of Sciences of Ukraine, 3b Naukova str., Lviv, 79000, Ukraine*

Correspondence should be addressed to Olexandr Polishchuk; od_polishchuk@ukr.net

**Abstract.** The methods are proposed for evaluation of complex dynamical systems, choice of their optimal operating modes, determination of optimal operating system from given class of equivalent systems, system's timeline behaviour analysis on the basis of versatile multicriteria and multilevel analysis of behaviour of system's elements.

## 1. Introduction

Study of complex dynamical systems of different types (technical, biological, social, economical) has been attracting the attention of many researchers for a long time already [1-4]. Important direction of such a study is development of methods for evaluation of state, operating quality and interaction between objects of those systems [5-10]. Among main problems arising in this case local and global evaluation of quality of complex dynamical systems [10, 11] may be pointed out, as well as determination of their optimal operating modes [12], analysis of system behaviour during certain period of time, the choice of optimal operating system from certain class of equivalent systems. We propose a unified approach to solving - problems listed above. It is based on their comprehensive (which presupposes consideration of as many characteristics of system's elements as possible), multicriteria and multiparameter analysis. In order to provide operative processing of results the developed evaluation is multilevel, which means formulation of conclusions of different generalization degree: from local ones that determine the behaviour of particular characteristics of system's elements to final ones that determine the quality of system's operation in general. Problems considered in this article and methods for their solution are illustrated with example of analysis of motion of man's musculoskeletal system (MMSS) with prosthetic lower limb [14-16].

## 2. Formulation of Problems

Let us consider dynamical system that consists of $N$ elements and is able to operate under $L$ modes. In order to simplify the explanation, let us assume that behaviour of every element of system under *l-th* mode is described by set of characteristics $A_{n,l,m}(t)$, $m = \overline{1,M}$, $l = \overline{1,L}$, where $n$ is the number of element, $n = \overline{1,N}$, $t \in [0, T]$, where $T$ is the duration of test research. Each of these characteristics is the result of experimental research or mathematical modelling of processes taking place in the system. To analyse the behaviour $A_{n,l,m}(t)$, let us use $K_l^m$ criteria. Let us denote with $\Omega_{n,l,m,k_l^m}$ and $\tilde{\Omega}_{n,l,m,k_l^m}$ the domains of reference and permissible values for $A_{n,l,m}(t)$ characteristic with regard to $k_l^m$ criterion, $\Omega_{n,l,m,k_l^m} \subset \tilde{\Omega}_{n,l,m,k_l^m}$, $k_l^m = \overline{1,K_l^m}$. The term "equivalent" implies systems with same content, type and destination, the law of motion for which on the sequence of operating modes is described by set of characteristics $A_{n,l,m}(t)$, $m = \overline{1,M}$, $l = \overline{1,L}$, $n = \overline{1,N}$, $t \in [0, T]$, that are to satisfy defined set of criteria. Let us denote with $G_Q$ class of equivalent dynamical systems which consists of $Q$ elements.

While modelling MMSS with prosthetic lower limb as the multilink system of solid bodies, its motion may be described through rhythmical, kinematical, dynamical, and energy characteristics etc [14-16]. Each of these characteristics is *n*-dimensional vector-function, components of which describe peculiarities of system's elements behaviour (separate joints of human body and prosthesis applied) in the course of motion. Among the evaluation criteria there may be pointed out the deviation from known average norm, walking asymmetry level, deviation from best current rehabilitation result etc. The domains of reference values for components of characteristics are represented by domains of their change in course of normal walking, 3% level of walking asymmetry, known data on best prosthesis results etc. Disabled persons of the same sex, age group and health condition with same level of lower limb amputation and same or different types of prostheses applied make up the class of equivalent systems.

In this work we shall consider following problems for complex systems evaluation.
1. *Evaluation of system's element operation quality.* Solution for this problem allows to determine elements

representing potential threat of general system operation failures and to analyse their impact on surrounding elements. For systems composed of elements of the same type, solution of this problem allows to determine elements operating in the best way, i.e. reference elements. Finally, development of generalized conclusions regarding general system operation quality is based on results of system's elements evaluation.

2. *Choice of optimal mode for system operation.* Solution of this problem allows determining both most "comfortable" and extreme system operating modes, as well as modes of potential failure.

3. *Evaluation of system operation quality.* Solution of this problem allows determining the general quality of system operation according to defined set of parameters, criteria and operating modes.

4. *Choice of optimally operating system from given class of equivalent systems.* Solution of this problem allows to determine the best (referential) or the worst systems of the class. Optimally operating elements, modes and systems determined in the evaluation process may be used along as practically reachable quality references.

5. *Analysis of system operation history.* Solution of this problem allows to track and forecast the quality of system operation, determine trends of its development in the context of improvement or deterioration and to prevent possible failures in advance.

Returning to the issue of human prosthesis analysis, results of elements operation quality evaluation allow defining MMSS joints that are exposed to overload in the course of motion and distort its kinematics. Choice of optimal mode allows determining most sparing pace of motion for disabled person. Evaluation of system operation quality may therefore be interpreted as the quality of prosthesis for specific invalid. Choice of optimal system from the class defines the best current prosthetic result which may be used as evaluation criterion. If the class is represented by the set "disabled person – set of prostheses of different constructions", the choice of optimal system means determination of most favourable prosthesis construction for specific patient. Analysis of system behaviour during certain period of time allows to evaluate the process of invalid adaptation to applied prosthesis etc. Study of evaluation results with the presence of negative or close to negative conclusions provides strong reasons for improvement or change in rehabilitation methods applied.

## 3. Parameters and Scales of Evaluation

We will evaluate behaviour of the characteristic $A_{n,l,m}(t)$ under $k_l^m$ -th criterion by means of parameters $h_{n,l}^{m,k,p} = \left\| \alpha_{n,l}^{m,k} \right\|_{H_p[0,T]}$, where $\alpha_{n,l}^{m,k}(t) = \rho(A_{n,l,m}(t), \Omega_{n,l,m,k_l^m})$ is a distance between $A_{n,l,m}(t)$ and domain of reference values for *n-th* component of this characteristic under criterion $k_l^m$. $H_p[0,T]$ is a line of functional spaces, e.g. $C_{p-1}[0,T]$, $W_2^{p-1}[0,T]$, $p = \overline{1, P_{n,l}^{m,k}}$, or their combination, $k_l^m = \overline{1, K_l^m}$, $m = \overline{1, M}$, $l = \overline{1, L}$, $n = \overline{1, N}$, $t \in [0, T]$. Parameter values in uniform metric allow to track separate peaks or disturbance in behaviour of provided characteristic and its derivatives, those ones in mean-squared metric allow to define average value of their falling beyond domains of reference or permissible values.

Local evaluation under the $h_{n,l}^{m,k,p}$ parameter is performed as follows. Let us denote with $h_{n,l,\min}^{m,k,p}$ and $h_{n,l,\max}^{m,k,p}$ its minimum and maximum permissible values of parameter accordingly. As usual, the value $h_{n,l,\min}^{m,k,p} = 0$ corresponds to characteristic $A_{n,l,m}(t) \in \Omega_{n,l,m,k_l^m}$ and the value $h_{n,l,\max}^{m,k,p} = \max\limits_{A_{n,l,m} \in \overline{\Omega}_{n,l,m,k_l^m}} h_{n,l}^{m,k,p}$. If continuous evaluation scale is accepted, value of local evaluation of behaviour for characteristic $A_{n,l,m}(t)$ under *k-th* criterion and *p-th* parameter is defined by means of formula

$$e_{C,n,l}^{m,k,p} = \upsilon(h_{n,l,\max}^{m,k,p} - h_{n,l}^{m,k,p})/(h_{n,l,\max}^{m,k,p} - h_{n,l,\min}^{m,k,p})$$

where $\upsilon$ is the normalizing coefficient, e.g. 10. Then the highest positive evaluation will correspond to characteristic with values not falling beyond domain of reference values, zero evaluation corresponds to characteristic with values reaching the limits of domain of permissible values, and negative evaluation corresponds to characteristic with values falling beyond domain of permissible values. Extent of this falling beyond limits is determined by absolute value of $e_{C,n,l}^{m,k,p}$.

If discrete evaluation scale is accepted, every real value of functional $h_{n,l}^{m,k,p}$ within the limits of interval $[h_{n,l,\min}^{m,k,p}, h_{n,l,\max}^{m,k,p}]$ corresponds to integer number. Let us suppose that $I$ is the number of grades of integer scale, e.g. 10 and $\delta_i \in [0,1], \delta_i < \delta_{i+1}, i = \overline{0(1)I}$, $\delta_0 = 0$, $\delta_{I+1} = 1$. Then integer rating evaluation for charac-teristic $A_{n,l,m}(t)$ is defined by coefficient $e_{D,n,l}^{m,k,p}=i$ if $e_{C,n,l}^{m,k,p}/\upsilon \in [\delta_i, \delta_{i+1}[, i=\overline{0(1)I}$.

If the number of grades of integer scale is not great (2–5), its values correspond to those ones of conceptual evaluation scale, where each grade of discrete scale corresponds to value "unsatisfactory", "satisfactory", "good", "excellent" in ascending order. It is obvious that under consecutive transition from continuous to conceptual scale, the values of evaluations become less distinctive. The latter one, as well as the discrete scale with low number of grades is practically unacceptable to trace insignificant changes in system's elements behaviour or forecasting of their behaviour. One more drawback of conceptual scale is the fact that grade "satisfactory" may denote any possible value from "almost good" to "slightly better than unsatisfactory". However, it is convenient and comprehensible in cases of single or rare routine system evaluations.

The hybrid scale seems to be most convenient for practical use. Hybrid scale is precise rating scale combining the benefits of continuous and conceptual scales. Let us develop it as follows. To evaluate characteristic $A_{n,l,m}(t)$ according to $k_l^m$ criterion by $n$-th component we will use the pair of parameters $h_{n,l}^{m,k,C} = \|\alpha_{n,l}^{m,k}\|_{C_0[0,T]}$ and $h_{n,l}^{m,k,L} = \|\alpha_{n,l}^{m,k}\|_{L_2[0,T]}$ that denote characteristic's deviation from domain of reference values in uniform and mean-squared metrics. Let's consider the characteristic evaluation $e_{D,n,l}^{m,k,C}$ "excellent" or equivalent to 5 if $h_{n,l}^{m,k,C}=0$, and "unsatisfactory" if $h_{n,l}^{m,k,C} >> h_{n,l,\max}^{m,k,C}$. Let us introduce such parameter $\delta \in (0,1)$ that integer rating evaluation $e_{D,n,l}^{m,k,C}$ under parameter $h_{n,l}^{m,k,C}$ is "satisfactory" or equivalent to 3 if continuous evaluation $e_{C,n,l}^{m,k,C}/2 \in (0,\delta)$, and "good" or equivalent to 4 if $e_{C,n,l}^{m,k,C}/2 \in [\delta, 1)$. Development of precise rating evaluations is illustrated by simple example. Let us consider that $\Omega_{n,l,m,k_l^m}=a_{ref}=const$, $t \in [0,T]$, $\tilde{\Omega}_{n,l,m,k_l^m}=[a_{ref}, a_{\max}] \times [0,T]$, $a_{\max}=const$, $\gamma=a_{ref}+\delta(a_{\max}-a_{ref})$, and $a(t)=\alpha_{n,l}^{m,k}(t)$, $t \in [0,T]$ (Fig. 1). Then we shall consider precise rating evaluations $e_{S,n,l}^{m,k,C}$ according to parameter $h_{n,l}^{m,k,C}$ that define presence and magnitude of disturbances in behaviour of characteristic $A_{n,l,m}(t)$ equivalent to

– 2, if $e_{D,n,l}^{m,k,C}=2$;

– 3 + ($a_{\max}-\|a(t)\|_{C_0[0,T]}$)/($a_{\max}-\gamma$), if $e_{D,n,l}^{m,k,C}=3$;

– 4 + ($\gamma-\|a(t)\|_{C_0[0,T]}$)/$\gamma$, if $e_{D,n,l}^{m,k,C}=4$;

– 5, if $e_{D,n,l}^{m,k,C}=5$.

Then we shall consider precise rating evaluations $e_{S,n,l}^{m,k,L}$ according to parameter $h_{n,l}^{m,k,L}$ that define mass character of disturbances in behaviour of characteristic $A_{n,l,m}(t)$ equivalent to

– 2, if $e_{D,n,l}^{m,k,C}=2$;

– 3 + (($a_{\max}-\gamma)\sqrt{T}-\|a(t)-\gamma\|_{L_2[0,T]}$)/($a_{\max}-\gamma)\sqrt{T}$, if $e_{D,n,l}^{m,k,C}=3$;

– 4 + ($\|\gamma-a(t)\|_{L_2[0,T]}/\gamma \sqrt{T}$, if $e_{D,n,l}^{m,k,C}=4$;

– 5, if $e_{D,n,l}^{m,k,C}=5$.

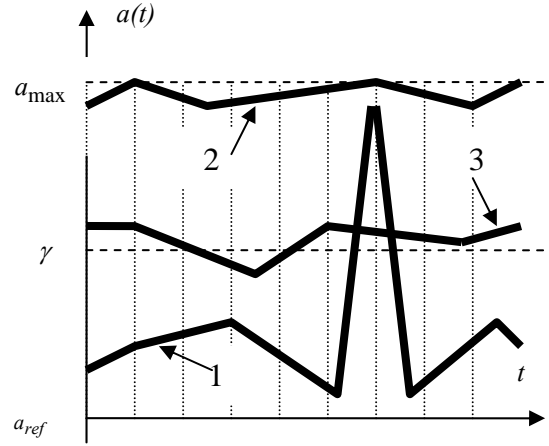

Fig.1

Then the pair of evaluations $e_{S,n,l}^{m,k,C}=3.05$, $e_{S,n,l}^{m,k,L}=3.98$ means the presence of low-number disturbances in behaviour of characteristic $A_{n,l,m}(t)$ (Fig.1, line 1). At the same time, the pair of evaluations $e_{S,n,l}^{m,k,C}=3.01$, $e_{S,n,l}^{m,k,L}=3.02$ shows that quality of element operation with regard to characteristic studied and evaluation criterion is close to critical (Fig.1, line 2). Pair of evaluations $e_{S,n,l}^{m,k,C}=3.95$, $e_{S,n,l}^{m,k,L}=3.91$ shows that quality of element operation with regard to characteristic studied and evaluation criterion is close to "good" (Fig.1, line 3). I. e. developed precise rating evaluati-ons provide quite specific, reasoned and understandable to average user information about the behaviour of evaluated characteristic of system element. Evaluations $e_{S,n,l}^{m,k,C}$ and $e_{S,n,l}^{m,k,L}$ hereinafter are referred to as "local". In the same way they may be developed for the whole range of functional spaces $H_p[0,T]$, $p = \overline{1, P_{n,l}^{m,k}}$, for arbitrary domains $\Omega_{n,l,m,k_l^m}$ and $\widetilde{\Omega}_{n,l,m,k_l^m}$.

In general, the number of numerical values of parameters of local evaluations of element is evaluated with the number $S_n = \sum_{l=1}^{L} \sum_{m=1}^{M} \sum_{k=1}^{K_l^m} P_{n,l}^{m,k}$, $n = \overline{1, N}$, and system elements selected for monitoring with the number $S = \sum_{n=1}^{N} S_n$. Let us take as an example motion of MMSS during walking [16]. The monitoring of behaviour of three pairs of MMSS joints (hip, knee, and ankle joints) is performed while they are implementing 18 functions (walking in slow, normal, and fast speed along horizontal or inclined surfaces (up and down) with and without load), each of which is described by three characteristics (kinematical, dynamical, and energy). When evaluating according to 4 criteria (deviation from known normal areas of human walking [16], level of motion asymmetry, deviation from best result reached, level of motion stability [15]) and 2 parameters (in uniform and mean-squared metrics), number of system elements evaluations $S_n = 432$, $n = \overline{1,6}$, general number of local evaluations $S=2592$.

## 4. Evaluation of system's element operation quality

It is obvious that direct analysis of totality of all numerical parameters of local evaluation is complex problem. For its solution, the sequence of weighted averaged evaluations of different generalization level is developed on the basis of local evaluations set until the final conclusion regarding the operation quality of considered system's element is reached.

Such development is performed for every element according to the whole set of parameters, criteria, characteristics and operating modes and presumes following levels of generalizations.

1) According to the set of parameters, for fixed criterion for evaluation of element's characteristic in operating mode specified:

$$H_{n,l}^{m,k} = (\rho_C\, e_{S,n,l}^{m,k,C} + \rho_L\, e_{S,n,l}^{m,k,L})/(\rho_C + \rho_L),$$
$$k = k_{n,l}^m = \overline{1, K_{n,l}^m},\; m = m_{n,l} = \overline{1, M_{n,l}},\; l = \overline{1, L},\; n = \overline{1, N},$$

where $\rho_C$, $\rho_L$ are weight coefficients defining evaluation parameters priority. Obtained value allows determining criteria under which evaluated element's characteristic in operating mode specified are unsatisfactory;

2) according to the set of evaluation criteria, for fixed element's characteristic in operating mode specified:

$$H_{n,l}^{m} = <\boldsymbol{\rho}^{Cr},\, \widetilde{\mathbf{H}}_{n,l}^{m}>_{R^{K_{n,l}^m}} / <\boldsymbol{\rho}^{Cr},\mathbf{1}>_{R^{K_{n,l}^m}},$$
$$m = m_{n,l} = \overline{1, M_{n,l}},\; l = \overline{1, L},\; n = \overline{1, N},$$

where $<.,.>_{R^K}$ is a scalar product in Euclidean space $R^K$, $\mathbf{1} = \{1\}_{k=1}^{K}$, $\widetilde{\mathbf{H}}_{n,l}^{m} = \{H_{n,l}^{m,k}\}_{k=1}^{K_{n,l}^m}$, $\boldsymbol{\rho}^{Cr} = \{\rho_k^{Cr}\}_{k=1}^{K_{n,l}^m}$ is vector of weight coefficients defining evaluation criteria priority. Obtained value allows determining characteristics under which operation of evaluated element in operating mode specified is unsatisfactory;

3) according to the set of element's characteristics in operating mode specified:

$$H_{n,l} = <\boldsymbol{\rho}^{Ch},\, \widetilde{\mathbf{H}}_{n,l}>_{R^{M_{n,l}}} / <\boldsymbol{\rho}^{Ch},\mathbf{1}>_{R^{M_{n,l}}},$$
$$l = \overline{1, L},\; n = \overline{1, N},$$

where $\widetilde{\mathbf{H}}_n^m = \{H_{n,l}^{m}\}_{m=1}^{M_{n,l}}$, $\boldsymbol{\rho}^{Ch} = \{\rho_m^{Ch}\}_{m=1}^{M_{n,l}}$ is vector of weight coefficients defining elements' characteristics priority. Obtained value allows determining operating modes in which operation of evaluated element is unsatisfactory;

4) for element specified according to set of operating modes:

$$H_n = <\boldsymbol{\rho}^{Mo},\, \widetilde{\mathbf{H}}_n>_{R^L} / <\boldsymbol{\rho}^{Mo},\mathbf{1}>_{R^L},\; n = \overline{1, N},$$

where $\widetilde{\mathbf{H}}_n = \{H_{n,l}\}_{l=1}^{L}$, $\boldsymbol{\rho}^{Mo} = \{\rho_l^{Mo}\}_{l=1}^{L}$ is vector of weight coefficients defining operating modes priority. Obtained value allows determining system elements, operation of which is unsatisfactory. Improvement and modification of those elements allows increasing general system operating quality. As for the problem of rehabilitation of disabled persons, unsatisfactory or close to unsatisfactory outcome of element evaluation means that selection of different prosthesis design parameters the is required. Those parameters should help to decrease load over retained joints of lower limbs and improve the kinematics of that person's motion.

## 5. Choice of optimal mode for system operation

If elements selected for monitoring possess have similar sets of evaluation parameters and criteria, as well as characteristics and operating modes, we are able to form generalized conclusions as to system's operation according to corresponding parameters, criteria, characteristics and operating modes. Let us develop the sequence of weighted average evaluations of different generalization degree on the basis of local evaluation set. This will allow analysing system's behaviour according to corresponding parameter, criterion or characteristic while operating in specified mode:

1) for separate evaluation parameter according to set of each characteristic's components for each evaluation criterion:

$$V_{l,C}^{m,k} = <\rho^{El}, e_{S,n,l}^{m,k,C}>_{R^N} / <\rho^{El}, \mathbf{1}>_{R^N},$$
$$V_{l,L}^{m,k} = <\rho^{El}, e_{S,n,l}^{m,k,L}>_{R^N} / <\rho^{El}, \mathbf{1}>_{R^N},$$
$$k = k_{n,l}^m = \overline{1, K_{n,l}^m}, \ m = m_{n,l} = \overline{1, M_{n,l}}, \ l = \overline{1, L},$$

where $\rho^{El} = \{\rho_n^{El}\}_{n=1}^N$ is vector of weight coefficient defining system's characteristics (elements) components priority;

2) for separate criterion according to set of evaluation parameters:

$$V_l^{m,k} = (\rho_C V_{l,C}^{m,k} + \rho_L V_{l,L}^{m,k})/(\rho_C + \rho_L),$$
$$k = k_{n,l}^m = \overline{1, K_{n,l}^m}, \ m = m_{n,l} = \overline{1, M_{n,l}}, \ l = \overline{1, L};$$

3) for separate characteristics according to set of evaluation criteria:

$$V_l^m = <\rho^{Cr}, \tilde{\mathbf{V}}_l^m>_{R^{K_{n,l}^m}} / <\rho^{Cr}, \mathbf{1}>_{R^{K_{n,l}^m}},$$
$$\tilde{\mathbf{V}}_l^m = \{V_l^{m,k}\}_{k=1}^{K_{n,l}^m}, \ m = m_{n,l} = \overline{1, M_{n,l}}, \ l = \overline{1, L};$$

4) for mode specified according to set of system elements characteristics:

$$V_l = <\rho^{Ch}, \tilde{\mathbf{V}}_l>_{R^{M_{n,l}}} / <\rho^{Ch}, \mathbf{1}>_{R^{M_{n,l}}},$$
$$\tilde{\mathbf{V}}_l = \{V_l^m\}_{m=1}^{M_{n,l}}, \ l = \overline{1, L}.$$

Evaluations $V_l$ allow determining operating modes in which system's operation is the worst. On the other hand, let us suppose that for system studied $K$ set of modes, $1 \leq K \leq l$, with the highest quality evaluation was received. Let us consider the procedure of optimal operating mode choice in case $K > 1$.

Choice of optimal system operating mode is made with consideration of following [17]. Let us presume that $\{a_l\}_{l=1}^L$ is optional set of real numbers under which $\sum_{l=1}^L a_l$. Among all numbers of this kind, value $\prod_{l=1}^L a_l$ reaches it's maximal point when $a_l = A/L$, $l = \overline{1, L}$. If $\{a_l\}_{l=1}^L$ is the set of evaluations, it means that deviation of their values from simple average is minimal.

Let us suppose that we have a few modes with highest evaluations. We shall consider mode (modes), for which maximum value of $\prod_{m=1}^{M_{n,l}} V_l^m$ is obtained, optimal operating modes for dynamical system studied. In the context of problem of rehabilitation of disabled persons, it means that efforts between the joints are distributed more evenly in a given mode of motion or kinematics of motion similar to the motion of the normal person.

## 6. Evaluation of system operation quality

Using previous outcomes, evaluation of system operating quality may be received in two ways. Thus, the value

$$H = <\rho^{El}, \tilde{\mathbf{H}}>_{R^N} / <\rho^{El}, \mathbf{1}>_{R^N}$$

where $\tilde{\mathbf{H}} = \{H_n\}_{n=1}^N$ provides global system evaluation, i. e. final conclusion regarding its operational quality.

The same evaluation will be received if the outcomes of system evaluation are generalized according to set of operating modes

$$V = <\rho^{Mo}, \tilde{\mathbf{V}}>_{R^L} / <\rho^{Mo}, \mathbf{1}>_{R^L}$$

where $\tilde{\mathbf{V}} = \{V_l\}_{l=1}^L$. It is obvious that $H \equiv V$.

## 7. Choice of system with optimal operation

Let us suppose that subclass $G_{\tilde{Q}}$, $1 \leq \tilde{Q} \leq Q$, of systems with highest evaluations of operating quality was obtained from $G_Q$ class. Let us consider the procedure of optimal system choice in case $1 < \tilde{Q}$, using algorithm applied for choice of optimal system operating mode. I. e., we will consider the system (systems) of set $G_{\tilde{Q}}$, for which maximum value of $\prod_{l=1}^L V_l$ is obtained, the optimal operating dynamical system of class. In context of problem of rehabilitation of disabled persons, this means that transition from one mode of motion to another causes the minimal inconvenience for invalid.

## 8. Analysis of change in system operation quality

Problem of analysis of change in systems operation quality in the timeline is similar to that one being considered in previous clause, the only difference being that class of equivalent systems is comprised by the same system but at different periods of time. Analysis for the best or the worst evaluation results will allow deter-mining most or least favourable conditions for its operation. If the sequence (prehistory) of system evaluations $\{V(T_j)\}_{j=1}^{J}, J \geq 2$, received at points of time $T_j$, $j = \overline{1, J}$, increases monotonically, the quality of system operation increases, if it decreases monotonically, the quality decreases, and if it is close to constant value $\sum_{j=1}^{J} V(T_j)/J$, the quality is stable. Evaluations history allows performing at least short-term forecasting of system operation quality. Indeed, let us suppose that $\mathbf{\Phi}(t) = \{\varphi_j(t)\}_{j=1}^{J}$ is the system of linearly independent functions, defined at the interval $[T_1, T_J]$. Let us develop function $V(t) = <\mathbf{A}, \mathbf{\Phi}(t)>_{R^J}$ where $\mathbf{A} = \{a_j\}_{j=1}^{J}$ is the vector of unknown coefficients. Then, forecasted value of system evaluation $V(t)$ at $T_{J+1} > T_J$ point of time is obtained from ratio $V(T_{J+1}) = <\mathbf{A}, \mathbf{\Phi}(T_{J+1})>_{R^J}$, where vector $\mathbf{A}$ is determined from condition $<\mathbf{A}, \mathbf{\Phi}(T_k)>_{R^J} = V(T_k)$, $k = \overline{1, J}$.

Prognostic analysis of precise evaluations allows determining point of time when conceptual evaluation will be reduced by one unit. In particular, when the sequence $\{V(T_j)\}_{j=1}^{J}$ is monotonically decreasing, the time for next system study may be defined from condition $V(t) \geq V^*$, where $V^*$ is the value corresponding to conceptual evaluation decreased by one unit comparing to that one determined at the moment of last examination.

## 9. Conclusions

In this work unified approach is proposed for evaluation of complex systems operations on all levels of their structuring. It is obvious, that substantiation of evaluation depends to the great extend, on integrity and completeness of both corresponding modes set as well as system elements characteristics, and set of evaluation criteria and parameters, as well as adequacy of weighted coefficients defining their priority. The great number of local and generalized evaluations requires development of convenient ways for evaluation outcomes visualization and disaggregation of global conclusions of different levels for localisation of reasons for drawbacks discovered [9, 14, 15].

Regarding the problems of disabled persons rehabilitation practices, methodology proposed allows, beside above described, to solve the problems of comparative analysis for age, sex and other peculiarities of human walking both in normal condition and with prosthetic lower limb, to study the impact of MMSS pathologies of different types on restriction of MMSS functional capabilities and impact of rehabilitation tools applied on reestablishment of such capabilities both in specific and in general cases, to perform comparative analysis of different rehabilitation methods etc.

## References

[1] A.-L. Barabási, "The architecture of complexity", *IEEE Control Systems Magazine*, vol. 27, no. 4, pp. 33-42, 2007.
[2] S. Boccatti, V. Latora, Y. Moreno, M. Chavez, and D.-U., "Complex Networks: Structure and Dynamics", *Physics Reports*, vol. 424, pp. 175-308, 2006.
[3] Y. Bar-Yam, "About Engineering Complex Systems: Multiscale Analysis and Evolutionary Engineering", *Engineering Self-Organising Systems: Methodologies and Applications*, Springer, London, England, pp. 16–31, 2005.
[4] D. Hinrichsen and A.J. Pritchard, *Mathematical Systems Theory*, Springer, New York, USA, 2005.
[5] J. Dombi, "Basic concepts for a theory of evaluation: the aggregative operator", *European Journal of Operational Research*, vol. 10, no. 3, pp. 282–293, 1982.
[6] A. Wittmuss, "Scalarizing multiobjective optimization problems", *Mathematical Researches*, vol. 27, pp. 255–258, 1985.
[7] C.L. Owen, "Evaluation of complex systems", *Designe Studies*, vol. 28, no. 1, pp. 73–101, 2007.
[8] D. Roy and T. Dasgupta, "Evaluation of reliability of complex systems by means of a discretizing approach", *International Journal of Quality & Reliability Management*, vol. 19, no. 6, pp. 792–801, 2002.
[9] L. Norros and P. Saviola, *Usability evaluation of complex systems*, STUK, Helsinki, Finland, 2004.
[10] D. Polishchuk, O. Polishchuk, and M. Yadzhak, "Comparison of methods of complex system evaluation", *Information Extraction and Processing*, vol. 32 (108), pp. 110-118, 2010.
[11] O. Polishchuk, "Optimization of function quality of complex dynamical systems", *Journal of Automation and Information Sciences*, no. 4, pp. 39-44, 2004.
[12] O. Polishchuk, "Choice of optimal regimes for functioning of complex dynamical systems", *Mathematical Methods and Physicomechanical Fields*, vol. 48, no. 3, pp. 62-67, 2005.
[13] O. Polishchuk, "About the choice of the optimal dynamical system from the class of equivalent systems", *Information Extraction and Processing*, vol. 20 (96), pp. 23-28, 2004.
[14] D. Polishchuk, O. Polishchuk, and M. Yadzhak, "Solution of some problems of evaluation of the complex systems", *Proceedings of the 15th International Conference on Automatic Control*, Odesa, 23-26 September 2008, pp. 968-976.
[15] O. Polishchuk, "Optimization of evaluation of man's musculo-sceletal system", *Computing Mathematics*, vol. 2, pp. 360-367, 2001.
[16] D.A. Winter, *The biomechanics and motor control of human gait: normal, elderly and pathological*, Univercity of Waterloo Press, Waterloo, Canada, 1991.
[17] G.H. Hardy, J.E. Littlewood, and G. Polia, *Inequalities,* Cambridge Univ. Press, Cambrige, England, 1988.